\def \version {2015--01--14}

\documentclass[12pt]{article}

\usepackage{amssymb}

\def \ex {\mbox{\rm ex}^*}
\def \cex {\overline{\mbox{\rm ex}^*}}

\def \bsk {\bigskip}
\def \nin {\noindent}

\def \pf {\nin{\bf Proof.} \ }
\def \qed {\hfill $\Box$}

\def \edik {^{\mbox{\scriptsize\rm th}}}
\def \ith {$i\edik$}

\def \cB {{\cal B}}
\def \cG {{\cal G}}

  \newtheorem{thm}{Theorem}
  \newtheorem{cor}[thm]{Corollary}
  \newtheorem{prp}[thm]{Proposition}
  
  \newtheorem{prm}[thm]{Problem}
  \newtheorem{cnj}[thm]{Conjecture}
  \newtheorem{lem}[thm]{Lemma}
  \newtheorem{dfn}[thm]{Definition}

\begin{document}

\title{Transversal designs and induced decompositions
 of graphs \thanks{~Research supported in part
   by the Hungarian State and the European Union
    under the grant TAMOP-4.2.2.A-11/1/KONV-2012-0072 } 
    }
\author{Csilla Bujt\'as~$^1$\qquad
   \vspace{2ex}
        Zsolt Tuza~$^{1,2}$\\
\normalsize $^1$~Department of Computer Science and Systems
 Technology \\
  \normalsize University of Pannonia \\
\normalsize  Veszpr\'em,
   \vspace{1ex}
     Hungary \\
\normalsize $^2$~Alfr\'ed R\'enyi Institute of Mathematics \\
       \normalsize Hungarian Academy of Sciences \\
\normalsize  Budapest,
   \vspace{1ex}
     Hungary
 }
\date{\small Latest update on \version}
\maketitle

\begin{abstract}
We prove that for every complete multipartite graph $F$
 there exist very dense graphs $G_n$ on $n$ vertices, namely with
 as many as ${n\choose 2} - cn$ edges for all $n$,
 for some constant $c=c(F)$,
 such that $G_n$ can be decomposed into edge-disjoint
 \emph{induced} subgraphs isomorphic to~$F$.
This result identifies and structurally explains
 a gap between the growth rates
 $O(n)$ and $\Omega(n^{3/2})$ on the minimum number of
 non-edges in graphs admitting an induced $F$-decomposition.

\bsk

 \noindent
\textbf{2010 Mathematics Subject Classification: }
 05C35, 
 05C70. 

\bsk

 \noindent
\textbf{Keywords and Phrases: } induced subgraph, edge decomposition,
 complete multipartite graph, extremal graph theory.

\end{abstract}

\vfill \newpage

\section{Introduction}

In this paper we study an extremal problem raised by Bondy and Szwarcfiter
\cite{BS} in 2013,  and extend the results proved in
\cite{BS,CT,HT}. Throughout, we consider simple undirected graphs,
 with the central concept of induced decomposition, under which we
 mean edge-decomposition into induced subgraphs.

 \subsection{The problem to be studied}
 \bsk
 \begin{dfn}\rm
  For two graphs\/ $F$ and\/ $G$, an \emph{induced\/ $F$-decomposition} of\/ $G$
 is a collection\/ $\{F_1,\dots , F_\ell\}$ of induced subgraphs of\/
 $G$ such that
  \begin{itemize}
  \item each\/ $F_i$ is isomorphic to\/ $F$\/ (\/$1\le i \le \ell$),
  \item the subgraphs are pairwise edge-disjoint, and
  \item $\bigcup_{i=1}^\ell E(F_i) = E(G)$.
   \end{itemize}
        \end{dfn}

\bsk

Clearly, for any fixed $F$ different from $K_2$, not all graphs $G$
admit an induced $F$-decomposition. Especially, if $F$ is not a
complete graph and the induced decomposition exists, then $G$ cannot
be complete. Therefore, the following problem was proposed by Bondy
and Szwarcfiter.
 \begin{prm}[\cite{BS}] Given a graph\/ $F$ and a positive integer\/ $n
 \ge |V(F)|$, determine the maximum number\/ $\ex(n,F)$ of edges in a graph\/ $G$
  of order\/ $n$ which admits an induced\/ $F$-decomposition.
 \end{prm}

 For the asymptotic version of this problem,  the following basic
 result was proved by Cohen and Tuza.

 \begin{thm}[\cite{CT}] \label{ex}
  For every graph\/ $F$ we have\/ $\ex(n,F)={n \choose 2}-o(n^2)$.
 \end{thm}

 Note that this  theorem is analogous to a result of
   Frankl and F\"uredi from~\cite{FF}; the latter one is discussed
   in a more general framework.

Theorem~\ref{ex}  suggests that in fact the `complementary' function
$$\cex(n,F) = {n\choose 2} - \ex(n,F)$$
is the appropriate subject   of study, as it is more informative
when asymptotic results are desired.

\subsection{A short summary of earlier asymptotic results}

\begin{itemize}
\item If $F$ is a complete graph  then, by Wilson's theorem \cite{W}, we
infer that for each $k$ the equality $\cex(n,K_k)=0$ holds with
infinitely many values of $n$. In general, $\cex(n, K_k)=O(n)$ holds
for every $k \ge 2$.
\item As observed in \cite{BS}, if $F'$ is obtained from $F$ by
extending it with an isolated vertex,
$$\cex(n,F) \le \cex(n,F') \le \cex(n-1,F)+n-1 $$
holds.
\item If $F$ is  isolate-free and non-complete  then $\cex(n,F)=\Omega(n)$.
This was proved recently by Hal\'asz and Tuza \cite{HT}.
Particularly, this linear lower bound is asymptotically sharp; that
is, $\cex(n,F)=\Theta(n)$ holds if $F$  is a complete bipartite or
tripartite  graph which is non-complete \cite{HT}.
\item If $F$ is isolate-free, non-complete and not a complete multipartite
graph, then $\cex(n,F)=\Omega(n^{3/2})$ \cite{HT}. This bound was
proved to be asymptotically sharp for $P_4$, $C_6$, $2K_2$ and
$K_4-P_3$ in \cite{BS,CT,HT}. Moreover, as a consequence of results
in \cite{AMS},   with Noga Alon we proved that
$\cex(n,F)=\Theta(n^{3/2})$ is valid for every disconnected graph
$F$ whose each component is complete bipartite \cite{AT}.
\end{itemize}
A more detailed summary of results (not only of the asymptotic ones)
 and of open problems is given in \cite{BT}.

\bsk
  Concerning the lower bound $\Omega(n)$ in the third point
above, one can observe
 that the exclusion of isolated vertices is in fact irrelevant.

\begin{prp}   \label{p:non-nbr}
 If each vertex of\/ $F$ has a non-neighbor (and, in particular, if\/
  $F$ has an isolated vertex), then\/ $\cex(n,F)\ge n/2$.
\end{prp}

\pf
 Under the given assumption, if $G$ admits an induced $F$-decomposition then
 every vertex of $G$, too, must have a non-neighbor.
\qed

\subsection{Our main contributions}



The most important contribution of this paper is as follows.

\begin{thm}   \label{fotetel}
 If\/ $F$ is a complete multipartite graph, then\/ $\cex(n,F)=O(n)$.
\end{thm}

Together with the earlier results cited above and with
 Proposition \ref{p:non-nbr}, Theorem~\ref{fotetel}
 immediately implies the following consequences.
The first one is the characterization of linear growth, while the
 second one identifies a gap in the exponents of possible
 growth functions for $\cex(n,F)$.

\begin{cor}
  We have\/ $\cex(n,F)=\Theta(n)$ if and only if either\/ $F$ is a
  complete graph plus at least one isolated vertex, or\/
  $F$ is a complete multipartite graph having at least one non-edge,
   possibly together with some isolated vertices.
\end{cor}

\begin{cor}
 If\/ $1 <c < 3/2$, there exists no\/ $F$ satisfying\/
$\cex(n,F)=\Theta(n^c)$.
\end{cor}

In some steps of the constructions we apply design-theoretic
 methods, in particular `transversal designs', whose formal
 definition is given in Section~\ref{sect2}.
Using those structures the following result will be proved.

\begin{thm}   \label{t:kpart}
 Let\/ $F = K_{a_1,\dots,a_k}$, and suppose that there exists a
  transversal design TD\/$(k,a_i)$ for each\/ $i=1,2,\dots,k$.
 Then the complete multipartite graph\/ $F^* = K_{ma_1,\dots,ma_k}$
  with\/ $m = \prod_{i=1}^k a_i$ admits an\/ $F$-decomposition.
 Moreover, such a decomposition can be encoded with\/ $m^2$ integer
  sequences of length\/ $2k$ whose\/ \ith\ and\/ $(k+i)\edik$
   terms vary between 1 and\/ $a_i$.
\end{thm}

After some preliminary observations in  Section~\ref{sect2}, we
 prove our main result, Theorem \ref{fotetel} in Section~\ref{sect3}.
  The proof of Theorem \ref{t:kpart} is given in Section~\ref{s:transv}.
We close the paper with several open problems and conjectures.




\section{Edge decompositions of complete $k$-partite graphs}
\label{sect2}

In this section we make some preparations towards the proofs
 of
  Theorems \ref{fotetel} and \ref{t:kpart}.
In the first subsection we recall some important known facts
 from the theory of transversal designs.
In the second subsection we construct decompositions of
 some restricted classes of graphs --- having the same number of
  partite classes that the fixed given graph $F$ has ---
  which will be applied later for the general results.

\subsection{Transversal designs}

%
%
%

In general, a
 \emph{transversal design} has three parameters:
 the ``order'' or ``groupsize'' $n$, the ``blocksize'' $k$, and
 the ``index'' $\lambda$.
For our purpose, however, only the case $\lambda = 1$ will be relevant,
 therefore we disregard the general definition.
Hence, the combinatorial structure denoted TD$(k, n)$
 is a triple $(X, \cG, \cB)$, where
 \begin{enumerate}
\item $X$ is a set of $kn$ elements;
\item $\cG = (X_1,X_2,\dots,X_k)$ is a partition of $X$ into $k$ classes $X_i$
 (the groups), each of size $n$;
\item $\cB$ is a collection of $k$-element subsets of $X$ (the blocks);
\item every unordered pair of elements from $X$ is contained
 either in exactly one group
 or in exactly one block,
  but not both.
 \end{enumerate}

\bsk

The more widely known structures of Latin squares are arrangements of the
 numbers $1,2,\dots,n$ in an $n\times n$ square in such a way that
 each number occurs precisely once in each row and also in each column.
Two Latin squares $(a_{ij})_{1\le i,j\le n}$ and $(b_{ij})_{1\le i,j\le n}$
 are said to be \emph{orthogonal} if their position-wise concatenation
 $\{a_{ij} b_{ij} \mid 1\le i\le n, \ 1\le j\le n\}$ lists all the $n^2$
 ordered pairs of $\{1,2,\dots,n\}\times \{1,2,\dots,n\}$.

%
%
%

We shall apply the following two facts, the first of which is
 not very hard to prove.

\begin{lem} (see Theorem 3.18 in \cite{CD})   \label{TD-Latin equiv}
 The existence of a TD\/$(k,n)$ is equivalent to the existence of\/
  $k-2$ mutually orthogonal Latin squares of order\/~$n$.
\end{lem}

\begin{lem} (Chowla, Erd\H os and Straus, \cite{CES})   \label{many MOLS}
 The maximum number of mutually orthogonal Latin squares of
  order\/ $n$ tends to infinity as\/ $n\to\infty$.
\end{lem}

\subsection{Decompositions of complete $k$-partite graphs}

Next, we describe some steps which make it possible to apply
 transversal designs to the general decomposition problem.
We begin with a definition.

\begin{dfn}\rm
 Let $G=(V,E)$ and $F=(X,H)$ be $k$-partite graphs with
  vertex partitions $V=V_1\cup\cdots\cup V_k$ and
  $X=X_1\cup\cdots\cup X_k$.
 Assume further that $|V|/|X|=p$ is an integer and that
  $|V_i|=p|X_i|$ for all $1\le i\le k$.
 An \emph{embedded\/ $F$-decomposition of\/ $G$} is an
  $F$-decomposition with the further property that
  each $V_i$ is partitioned into $p$ sets,
  $$
    V_i = V_{i,1} \cup \cdots \cup V_{i,p}
  $$
  and, for each copy of $F$ in the decomposition,
  each partite set $X_i$
  coincides with some $V_{i,j}$ inside  $V_i$.
\end{dfn}

\begin{lem}   \label{l:embed}
 Let\/ $k$ and\/ $a_1,\dots,a_k$ be fixed natural numbers.
 Then, for every sufficiently large integer\/ $p$,
  the complete multipartite graph\/ $K_{pa_1,\dots,pa_k}$
  has an embedded\/ $K_{a_1,\dots,a_k}$-decomposition.
\end{lem}

 \pf
Based on Lemma \ref{many MOLS} we can choose a threshold value $p_0=p_0(k)$
 such that there exist at least $k-2$ mutually orthogonal Latin squares
 of order $p$ whenever $p>p_0$ holds.
Then, by Lemma \ref{TD-Latin equiv}, there exists a transversal design
 TD$(k,p)$ of order $p$ and blocksize $k$.
Observe that TD$(k,p)$
 is precisely an edge decomposition of the complete $k$-partite
 graph $K_{p,\dots,p}$ into copies of the complete graph $K_k$.
Substituting mutually disjoint sets of cardinality $a_i$ for
 all vertices in the $i\edik$ group of TD$(k,p)$ for $i=1,2,\dots,k$,
  each copy of $K_k$ becomes a copy of $K_{a_1,\dots,a_k}$.
Thus, an embedded\/ $K_{a_1,\dots,a_k}$-decomposition of
 $K_{pa_1,\dots,pa_k}$ is obtained.
 \qed

\bsk

Choosing $p_0$ as in the proof above, and letting $a^*_i = pa_i$
 for $i=1,2,\dots,k$ with an arbitrarily fixed $p>p_0$, we obtain:

\begin{cor}   \label{large multipl}
 For every\/ $k$-tuple\/ $(a_1,\dots,a_k)$ of positive integers
  there exists a\/ $k$-tuple\/ $(a^*_1,\dots,a^*_k)$ such that
  \begin{itemize}
   \item $K_{a^*_1,\dots,a^*_k}$ has an
    embedded\/ $K_{a_1,\dots,a_k}$-decomposition;
   \item a transversal design TD\/$^i =$\,TD\/$(k,a^*_i)$ exists
    for each\/ $i=1,2,\dots,k$.
  \end{itemize}
\end{cor}

The applicability of this approach is demonstrated by Theorem \ref{t:kpart},
 which was stated in the Introduction and will be proved in Section \ref{s:transv}.

\section{Proof of the linear upper bound} \label{sect3}

In this section we prove Theorem \ref{fotetel}.

 Let $F = K_{a_1,\dots,a_k}$ be a fixed given graph.
With reference to Lemma \ref{l:embed} we choose an integer $p>1$
 such that
  \begin{itemize}
   \item $p$ is a multiple of $\prod_{i=1}^k a_i$,
   \item $F_p:=K_{pa_1,\dots,pa_k}$ has an embedded $F$-decomposition.
  \end{itemize}

The general form of Wilson's theorem yields that there exists a
 modulus $q$ and a residue $r$ such that the complete graph $K_{sq+r}$
 has an $F$-decomposition whenever $s\ge s_0=s_0(F)$ is sufficiently large.
Note that, regarding the linear upper bound to be proved, we may
 assume that $n$ is large.
We write $n$ in the form
 $$
   n = (sq+r)\cdot p + t \qquad (s\ge s_0, \ \ 0\le t\le pq-1)
 $$
 and define
 $$
   n' = sq+r .
 $$

\paragraph{Step 1.}

 Decompose $K_{n'}$ into edge-disjoint copies of $F$.
Of course, those copies are non-induced subgraphs of $K_{n'}$
 (unless $F$ itself is a complete graph, which case we disregard).

\paragraph{Step 2.}

Replace each vertex of $K_{n'}$ with an independent set
 of cardinality $p$.
That is, two vertices are adjacent if and only if they belong to
 different independent sets;
 in this way the graph
  $$
    G' := K_{(sq+r)\cdot p} - (sq+r)\cdot K_p
      = K_{n-t} - \textstyle \frac{n-t}{p}\cdot K_p
  $$
 is obtained.
Moreover, each copy of $F$ in $K_{n'}$ becomes a copy of $F_p$
 in $G'$.
Observe that, although the copies of $F_p$ are still non-induced,
 the \ith\ vertex class of cardinality $pa_i$ is the union of
 $a_i$ independent sets of size $p$ each.

\paragraph{Step 3.}

In any copy of $F_p$, if an independent set $S$ of size $p$ belongs
 to the \ith\ vertex class of $F_p$, then partition $S$ into
 $p/\!a_i$ disjoint sets.
In this way the $a_i$ independent $p$-element sets inside the \ith\ vertex class
 together yield a partition of that class into $p$
 independent sets of size $a_i$ each.

\paragraph{Step 4.}

Denoting by $V_1,\dots,V_k$ the vertex classes in the copy of $F_p$,
 and by $V_{i,j}$ the sets in the partition obtained
 ($1\le i\le k$, $1\le j\le p$, $|V_{i,j}|=a_i$ for all $i$ and $j$),
 the embedded $F$-decomposition of $F_p$
  --- guaranteed by Lemma \ref{l:embed} and by the choice of $p$ ---
 yields $p^2$ copies of $F$ which are \emph{induced} subgraphs of $G'$.

\paragraph{Step 5.}

Complete the construction with $t$ isolated vertices,
 to obtain a graph $G$.
It is clear that this $G$ admits an induced $F$-decomposition.

\bsk

It remains to observe that the complement of $G$ has no more
 than a linear number of edges.
There are ${t\choose 2} + tn - t^2 < tn < pqn$
 vertex pairs containing at least one of the $t$ isolated vertices;
 and the number of non-edges in $G'$ is exactly
 $\frac{n-t}{p}\cdot{p\choose 2} < pn/2$.
Thus,
 $$
   \cex(n,F) \le {n\choose 2} - |E(G)| < (pq + p/2)\cdot n = O(n) .
 $$
This completes the proof of Theorem \ref{fotetel}.

\section{Complete $k$-partite graphs and transversal designs}
   \label{s:transv}

 Here we prove Theorem \ref{t:kpart}.
Recall that $F = K_{a_1,\dots,a_k}$, $F^* = K_{ma_1,\dots,ma_k}$,
 and $m = \prod_{i=1}^k a_i$.
It is now assumed that there exists a transversal design
  TD$(k,a_i)=(X^i,\cG^i,\cB^i)$ for every $i=1,2,\dots,k$,
   which we shall refer to as TD$^i$ or TD$^i(k,a_i)$.
Those TD$^i$ will become important in the last step of the
 construction.

In order to simplify notation throughout,
 we write the usual shorthand $[a_i]=\{1,2,\dots,a_i\}$.

\bsk

 \nin
 {\bf Notation.} \
For $1\le i\le k$, $1\le p\le a_i$, and $1\le q\le k$,
 the $p\edik$ element in the $q\edik$ group $X_q^i$ of TD$^i(k,a_i)$
 will be denoted by $\pi^i(p,q)$.
  (Inside each group the elements are numbered as $1,2,\dots,a_i$.)
  Note that for all $1\le p,p'\le a_i$ and all $1\le q,q'\le k$ ($q\ne
q'$),
 a block in $\cB^i$ containing the $p\edik$ element
 of the $q\edik$ group and the $(p')\edik$ element of the
 $(q')\edik$ group in TD$^i(k,a_i)$
 \emph{exists and is unique}.

\paragraph{Construction.}

Let the vertex set of $F^*$ be
$$
  V = V^1 \cup\cdots\cup V^k,
    \qquad |V^i| = ma_i \ \ \mbox{\rm for} \ \ i = 1,\dots,k.
$$
Perhaps the best way of thinking about $V^i$ is to view it as a
$k$-dimensional box
 composed of unit-cube cells containing $a_i$ items each; in the $j\edik$ direction
 the box has size $a_j$.
In this view we see $a_i$ layers orthogonal to the $i\edik$ direction, each of those
 layers consists of $m/a_i$ cells, and hence each layer contains precisely $m$ items.

Stating this in a more formal description,
 we partition each $V^i$ into $m$ disjoint sets of cardinality $a_i$
 (corresponding to the contents of the cells), and
 represent each of those subsets inside $V^i$ with a vector
$$
  {\mathbf j} = (j_1,\dots,j_k) \in [a_1] \times \cdots \times [a_k]
$$
 of length $k$, so that
$$
  V^i = \bigcup_{{\mathbf j} \in [a_1] \times \cdots \times [a_k]}
    V^i_{\mathbf j} .
$$
In this way a copy of $F$ inside $F^*$ is characterized by a
 $k$-tuple of vectors,
$$
  ({\mathbf j^1},\dots,{\mathbf j^k}) \in
    ([a_1] \times\cdots\times [a_k])^k ,
$$
 where
$$
  {\mathbf j^i} = (j_1^i,\dots,j_k^i) \in [a_1] \times \cdots \times [a_k]
$$
 for $i=1,2,\dots,k$.
Let us refer to this $({\mathbf j^1},\dots,{\mathbf j^k})$ as the
 \emph{detailed representation} of the copy in question.

There are $(ma_i)(ma_{i'})=m^2a_ia_{i'}$ edges of $F^*$ between $V^i$
 and $V^{i'}$, and each copy of $F$ contains $a_ia_{i'}$ edges between
 those two partite sets.
  Thus, we have to define $m^2$ edge-disjoint copies of $F$ inside $F^*$.
The key point is how to shuffle the vertex sets of the $m^2$ copies to ensure
 that those subgraphs compose a partition of the edge set of $F^*$.

The copies of $F$ will be encoded with vectors of length $2k$
 (i.e., a suitable collection of $k^2$-dimensional vectors
 will be encoded with $(2k)$-dimen\-sional ones),
 which we shall view as concatenations of two vectors
 $$
   ({\mathbf b}, {\mathbf c}) \in ( [a_1] \times \cdots \times [a_k] )^2
 $$
where
$$
  {\mathbf b} = (b_1,\dots,b_k) \in [a_1] \times \cdots \times [a_k],
$$
$$
  {\mathbf c} = (c_1,\dots,c_k) \in [a_1] \times \cdots \times [a_k]
$$
 are vectors of length $k$.

We identify a copy of $F$ for any given
 $({\mathbf b}, {\mathbf c})$ in the way detailed next.
\begin{itemize}
 \item
  Generally speaking, the
   first $k$ coordinates (i.e., those of ${\mathbf b}$) specify
  a $k$-tuple of layers (one layer in each of the boxes $V^i$) that means
  a copy of $K_{m,\dots,m}$; and the role of ${\mathbf c}$ is similar,
  but ``rotating'' by 1 among the boxes, e.g.\ while the first coordinate
  of ${\mathbf b}$ will provide a data for $V^1$, the first coordinate
   of ${\mathbf c}$ will be relevant for $V^2$, hence specifying a
   second copy of $K_{m,\dots,m}$.
  After that, the transversal designs will be used for defining further
   $k-2$ layers inside each $V^i$.
  Having all this at hand, exactly $k$ layers are determined inside
   each box, and all the intersections of those layers
    in the different boxes together uniquely identify
    a copy of $F$ associated with $({\mathbf b},{\mathbf c}$).
 \item The $i\edik$ coordinate $b_i$ of ${\mathbf b}$ determines that
  the $i\edik$ partite set of $F$ (inside $V^i$)
  must be one of those
  $V^i_{{\mathbf j}}$ which have $j_i=b_i$.
  This equivalently means that $j^i_i=b_i$ holds in the detailed representation.
 (Recall that ${\mathbf j^i}$ has the form
   ${\mathbf j^i}=(j_1^i,\dots,j_k^i)$.)
%
 \item The $i\edik$ coordinate $c_i$ of ${\mathbf c}$ determines the
  $i\edik$ coordinate $j_i^{i+1}$ of the next ${\mathbf j^{i+1}}$
   in the detailed representation.
  (For $i=k$ we mean ${\mathbf j^{k+1}}={\mathbf j^{1}}$.)
%
%
 \item By the definition of the transversal designs TD$^i(k,a_i)$, there is
   a unique block $B^i=B^i({\mathbf b}, {\mathbf c})\in\cB^i$ such that
   both $\pi^i(b_i,i)\in B^i$ and $\pi^i(c_i,{i+1})\in B^i$.
  (Also here, the successor $i+1$ of $i=k$ is meant to be 1.)
  Then the $i\edik$ coordinate of the copy of $F$ associated with
   $({\mathbf b}, {\mathbf c})$ in $V^{i'}$, that is $j^{i'}_i$
    in the representation $({\mathbf j^1},\dots, {\mathbf j^k})$,
    is defined as $j^{i'}_i=B^i\cap X_{i'}^i$ for all $i'$.
  This is our defining rule for $i'\notin\{i,i+1\}$, but the formula is
   valid also for $i'\in\{i,i+1\}$ because this was the way
   we selected $B^i$.
\end{itemize}

\paragraph{Verification of requirements.}

 It is clear that precisely $m^2$ copies are defined in this way.
Since the total number of edges in those copies is equal to that
 in $F^*$, the proof will be done if we show that each edge of $F^*$
 is contained in some copy of $F$.
This is equivalent to saying that each edge $v_iv_{i'}$
 ($v_i\in V^i$, $v_{i'}\in V^{i'}$) determines all the $k^2$ coordinates
  of some $k$-tuple of vectors $({\mathbf j^1},\dots, {\mathbf j^k})$.
The property is valid indeed, as can be seen from the following
 considerations.

The positions of $v_i$ and $v_{i'}$ determine the representing vectors
 ${\mathbf j^i}$ and ${\mathbf j^{i'}}$ of the $a_i$-element and
 $a_{i'}$-element sets
 (cells) containing $v_i$ and $v_{i'}$, respectively.
This determines $k$ pairs of coordinates:
$$
  (j_1^i,j_1^{i'}), \quad(j_2^i,j_2^{i'}), \quad \dots, \quad(j_k^i,j_k^{i'}) .
$$
The pair $(j_\ell^i,j_\ell^{i'})$ uniquely determines the block
 $B^\ell=B^\ell({\mathbf b}, {\mathbf c})\in\cB^\ell$
 in TD$^\ell(k,a_\ell)$ for $\ell=1,\dots,k$;
 this $B^\ell$ is the unique block containing the $(j_\ell^i)\edik$ element
  of the $i\edik$ block and the $(j_\ell^{i'})\edik$ element
  of the $(i')\edik$ block in TD$^\ell(k,a_\ell)$.
Thus, all the $k^2$ coordinates of the vectors
 ${\mathbf j^1},\dots,{\mathbf j^k}$ can be computed
 from those $k$ blocks by the rule $j_\ell^i=B^\ell\cap X^\ell_i$.

 This completes the proof of the theorem.

\section{Concluding remarks}

In this paper we studied the growth of the number of edges in
 graphs which admit edge decompositions into induced subgraphs
 isomorphic to a given complete multipartite graph.
  We proved a linear upper bound on the number of non-edges.
More explicitly, our construction yields the following.

\begin{thm}
 If\/ $n$ satisfies the following three conditions:
\begin{itemize}
 \item ${n\choose 2}$ is divisible by the number\/
   $\displaystyle\sum_{1\le i<i'\le k} a_ia_{i'}$ \
    of edges in $F$,
 \item $n-1$ is divisible by the greatest common divisor
  of the vertex degrees\/
   \ $\left( \sum_{i=1}^k a_i \right) -a_1$, \ $\left( \sum_{i=1}^k a_i \right) -a_2$,
   \ $\dots,$ \ $\left( \sum_{i=1}^k a_i \right) -a_k$,

 \item $n>n_0$ for some threshold value\/ $n_0=n_0(a_1,\dots,a_k)$,
\end{itemize}
 then the complete\/ $n$-partite graph\/ $K_{m,\dots,m}$
 with\/ $m=\prod_{i=1}^k a_i$ admits an
 edge decomposition into induced subgraphs isomorphic to\/ $F$.
\end{thm}

Comparison with the results from \cite{HT}
 identifies a gap in the exponent of $n$;
 this can be expressed in the following way, where
 the first case is trivial and
 the second case is taken from Wilson's theorem.

\vbox{
\emph{
\begin{itemize}
 \item[] \hspace{-1.2em}Let\/ $F$ be a graph without isolated vertices.
 \item If\/ $F=K_2$, then\/ $\cex(n,F)=0$.
 \item If\/ $F$ is a complete graph of order at least 3, then\/ $\cex(n,F)$
  oscillates between\/ $0$ and\/ $\Theta(n)$, reaching both extremities
  on sets of positive upper density.
 \item If\/ $F$ is a complete multipartite graph, but not a
  complete graph, then\/ $\cex(n,F)=\Theta(n)$.
 \item If\/ $F$ is not a complete multipartite graph, then\/
  $\cex(n,F)=\Omega(n^{3/2})$.
\end{itemize}
 }
 }

The next task could be to identify further gaps in the exponent of $n$
 describing the growth of $\cex(n,F)$, if there are any.
We believe that not all reals (and not all rational numbers)
 between $3/2$ and 2
 can occur as exponents for some $F$,
 and also that the gap between $n$
 and $n^{3/2}$ is not the only one.
This is expressed in the following sequence of conjectures.

\begin{cnj}
 There exist graphs\/ $F$ for which
  $$\lim_{n\to\infty} \frac{\cex(n,F)}{n^{3/2}} = \infty .$$
\end{cnj}

\begin{cnj}   \label{expon}
 For every non-complete graph\/ $F$, without isolated vertices,
  there exists a constant\/ $c=c(F)$ such that
  $\cex(n,F) = \Theta(n^c)$.
\end{cnj}

\begin{cnj}
 The constant\/ $c=c(F)$ in Conjecture \ref{expon} is rational
  for every graph\/ $F$.
\end{cnj}

\begin{cnj}   \label{c:nagys}
 Every non-complete graph\/ $F$ without isolated vertices satisfies
  the asymptotic equality \ $\cex(n,F)=\Theta(n^{2-1/t})$
  where\/ $t=t(F)$ is a positive integer.
\end{cnj}

If Conjecture~\ref{c:nagys} is true, the next question
is: Which structural property of $F$ determines $t(F)$\,?

\begin{cnj}   \label{c:aszimp}
 Every non-complete graph\/ $F$ without isolated vertices,
  the quotient
   $$\frac{\cex(n,F)}{n^{2-1/t}}$$
  tends to a positive limit as\/ $n$ gets large,
  for some positive integer\/ $t=t(F)$.
\end{cnj}

If this is true, we would be interested in the exact value of
the limit.

\begin{prm}
 Verify Conjectures \ref{c:nagys} and \ref{c:aszimp} for all graphs
  on at most six vertices.
\end{prm}

The smallest open case for Conjecture \ref{c:nagys} is $F=C_5$,
 while the one for Conjecture \ref{c:aszimp} is $F=2K_2$.

\paragraph{Hypergraphs.}

It is very natural to extend the above questions to
 $r$-uniform hypergraphs ($r\ge 3$), but rather too little is known
 so far about the analogous function $\cex(n,{\cal F})$.
Already the study of some very small examples may provide
 interesting pieces of information.

\end{document}